\documentclass{amsart}

\usepackage{amsmath,amsfonts,amssymb,amsthm,hyperref}
\usepackage{color}
\usepackage{centernot}

\usepackage[dvips]{graphicx}
\usepackage{latexsym}

\newtheorem{thm}{Theorem}[section]
\newtheorem{prop}[thm]{Proposition}

\newtheorem{lemma}[thm]{Lemma}

\newtheorem{cor}[thm]{Corollary}

\newtheorem{definitiontemp}[thm]{Definition}
\newenvironment{defn}{\begin{definitiontemp}
\normalfont}{\end{definitiontemp}}

\newcommand\aut[1]{\text{Aut}(#1)}

\def\bec{\begin{cor}}
\def\enc{\end{cor}}

\def\bet{\begin{thm}}
\def\ent{\end{thm}}

\def\becor{\begin{cor}}
\def\encor{\end{cor}}

\def\bel{\begin{lem}}
\def\enl{\end{lem}}

\def\bedef{\begin{defn}}
\def\endef{\end{defn}}

\def\bep{\begin{prop}}
\def\enp{\end{prop}}

\newenvironment{pf}{\begin{trivlist}\item[\hskip\labelsep
{\it Proof.}]}{\end{trivlist}}

\newcommand{\set}[2]{\ensuremath{ \{ #1 : #2 \} }}

\renewcommand{\deg}[1]{\ensuremath{\text{deg}(#1)}}

\newcommand{\F}{\mathfrak{F}}
\newcommand{\G}{\mathfrak{G}}
\newcommand{\Q}{\mathbb{Q}}
\newcommand{\R}{\mathbb{R}}
\newcommand{\C}{\mathbb{C}}
\newcommand{\A}{\mathcal{A}}

\newcommand{\U}{\mathcal{U}}

\newcommand{\Qbar}{\overline{\mathbb{Q}}}

\newcommand{\Avec}{\vec{A}}
\newcommand{\avec}{\vec{a}}

\newcommand{\fvec}{\vec{f}}
\newcommand{\Fvec}{\vec{F}}
\newcommand{\Gvec}{\vec{G}}
\newcommand{\gvec}{\vec{g}}
\newcommand{\Hvec}{\vec{H}}
\newcommand{\hvec}{\vec{h}}

\newcommand{\xvec}{\vec{x}}

\newcommand{\etavec}{\vec{\eta}}
\newcommand{\gammavec}{\vec{\gamma}}
\newcommand{\rhovec}{\vec{\rho}}

\newcommand{\TQbar}{T_{\Qbar}}
\newcommand{\ab}{\aut{\Qbar}}
\newcommand{\abz}{{\text{Aut}_{\bfz}(\Qbar)}}
\newcommand{\abb}{{\text{Aut}_{\bfb}(\Qbar)}}
\newcommand{\abc}{{\text{Aut}_{\bfc}(\Qbar)}}
\newcommand{\abd}{{\text{Aut}_{\bfd}(\Qbar)}}
\newcommand{\abI}{{\text{Aut}_{I}(\Qbar)}}
\newcommand{\abIA}{{\text{Aut}_{I_A}(\Qbar)}}
\newcommand{\an}{\aut{F_n}}
\newcommand{\am}{\aut{F_m}}
\newcommand{\as}{\aut{F_s}}

\newcommand{\Gal}[2]{\text{Gal}(#1/#2)}
\newcommand{\la}{\langle}
\newcommand{\ra}{\rangle}

\newcommand{\WKL}{\ensuremath{\textbf{WKL}_0}}
\newcommand{\PA}{\textbf{PA}}

\def\converges{\!\downarrow}

\renewcommand{\qed}{\hbox to 0pt{}\nobreak\hfill\rule{2mm}{2mm}}

\newcommand{\bfa}{\boldsymbol{a}}
\newcommand{\bfb}{\boldsymbol{b}}
\newcommand{\bfc}{\boldsymbol{c}}
\newcommand{\bfd}{\boldsymbol{d}}
\def\bfz{\boldsymbol{0}}
\def\s01{\ensuremath{\Sigma^0_1}}
\def\d02{\ensuremath{\Delta^0_2}}
\def\phi{\varphi}

\def\res{\!\!\upharpoonright\!\!}

\newcommand{\comment}[1]{}

\begin{document}

\title[Computability for the absolute Galois group of $\Q$]{Computability for the\\absolute Galois group of $\Q$}

\author{Russell Miller}
\thanks{
This work was initiated during a program supported
by the National Science Foundation under Grant \# DMS-1928930 and hosted by the
Simons Laufer Mathematical Sciences Institute in Berkeley, California, during the summer of 2022.
The author was partially supported by Grant \#581896 from the Simons Foundation,
and by several grants from The City University of New York PSC-CUNY Research Award Program.
}
\address{
\parbox{\dimexpr\textwidth -\parindent}{ \hspace*{0.1in}  Department of Mathematics, Queens College -- City University of New York, 65-30 Kissena Blvd., Flushing, NY 11367 USA \\
\hspace*{0.1in}  Ph.D. Programs in Mathematics and Computer Science, Graduate Center - City University of New York, 365 Fifth Avenue, New York, NY 10016 USA}}
\email{Russell.Miller@qc.cuny.edu}
\urladdr{\url{https://qcpages.qc.cuny.edu/~rmiller/}}

\begin{abstract}
The absolute Galois group $\Gal{\Qbar}{\Q}$ of the field $\Q$ of rational numbers can be
presented as a highly computable object, under the notion of type-2 Turing computation.
We formalize such a presentation and use it to address several effectiveness
questions about $\Gal{\Qbar}{\Q}$:  the difficulty of computing Skolem
functions for this group, the arithmetical complexity of various definable subsets
of the group, and the extent to which countable subgroups defined by complexity
(such as the group of all computable automorphisms of the algebraic closure $\Qbar$)
may be elementary subgroups of the overall group.
\end{abstract}

   
\maketitle

\section{Introduction}
\label{sec:intro}

The absolute Galois group $\Gal{\Qbar}{\Q}$ of the field $\Q$ of rational numbers
is traditionally and naturally viewed as a profinite group.  It is the object of intense study
in number theory and is often viewed as a dense and impenetrable object,
not well understood as yet.
In \cite{S94} it is stated that it is impossible to ``write down'' an element
of $\Gal{\Qbar}{\Q}$ distinct from the identity and complex conjugation.
Indeed, there is sometimes debate over the extent to which one can even
consider the algebraic closure $\Qbar$ itself in a reasonable way.

Nevertheless, from the point of view of computability theory, there is
a perfectly nice presentation of $\Qbar$ -- unique up to computable isomorphism --
from which one then can describe a very reasonable presentation of its
automorphism group $\ab$, which is to say, of $\Gal{\Qbar}{\Q}$.
This automorphism group is not susceptible to standard
computable structure theory, which focuses on the presentation of
countable structures.  However, in Section \ref{sec:presentation} we will present it
as the set of paths through a computable, finite-branching tree of height $\omega$.
Among all presentations of uncountable structures commonly occurring in mathematics,
our presentation of $\ab$ is probably the most congenial and tractable currently known,
far easier to work with than the usual presentations of the fields
$\R$ and $\C$, for instance.  Even among automorphism groups
of countable structures, it stands out as a particularly straightforward example.
Of course, it is for the number theorists and others studying $\ab$
for their own purposes to decide whether this presentation is of any use
to them or not.  Here we simply regard $\ab$ as a first-order structure
in the language of groups, uncountable but with a highly accessible
tree presentation, to be studied from the standpoint of computability
in the same manner as any other such structure would be.

We will see that there is a natural notion of the Turing degree of each automorphism
in $\ab$:  some are computable, most are not, but each one
has a well-defined spot in the usual pantheon of Turing degrees.
Indeed, for each Turing degree $\bfd$, the set of $\bfd$-computable automorphisms
$$ \abd = \set{f\in\ab}{f\leq_T\bfd}$$
forms a subgroup of $\ab$.  Moreover, for each number field $F$,
every element of $\Gal{F}{\Q}$ is the restriction of an element of $\abz$
(and likewise with every $\abd$).  Subgroups such as $\abz$ therefore
feel like very reasonable ``simulations'' of the full group $\ab$, which is the inverse
limit of the directed system of all such finite groups $\Gal{F}{\Q}$.  It is natural
to ask to what extent each $\abd$ is an elementary subgroup of the full group
$\ab$, or at least to what extent they share the same theory
in the language of groups.

Of course, it is not difficult to prove the existence 
of countable, fully elementary subgroups of $\ab$.
The Downward L\"owenheim-Skolem Theorem gives examples immediately.  
Our long-term goal is not so much to exhibit countable elementary subgroups
as to examine the extent to which naturally defined subgroups
such as these $\abd$ are elementary.  Another set of examples,
to which our investigations will lead us, are the unions $\cup_{\bfd\in I}\abd$
over various Scott ideals $I$ in the Turing degrees.  These prove more fruitful than
the individual groups $\abd$, which correspond to principal Turing ideals.
(For the relevant definitions, see Subsection \ref{subsec:lowbasis} of the Appendix.)
In Section \ref{sec:elementarity} we will see that the Scott-ideal groups
are elementary within $\ab$ for all universal formulas, all existential formulas,
all positive formulas, and a significant class of more complex formulas.

Our investigations of elementarity will place us in a good position,
in Section \ref{sec:definability}, to consider the complexity of subsets
of $\ab$ definable by formulas of the types listed above.
Some of these sets will prove to have complexity significantly lower than
the complexity of the defining formula.  Additionally, in Section \ref{sec:Skolem},
we will apply the techniques developed for elementarity to consider the complexity
of Skolem functions for $\ab$:  when $\alpha(F,G)$ is a formula in the language
of groups and it is known that $(\forall F\exists G)\alpha(F,G)$ holds in $\ab$,
how difficult is it, given a specific $f\in\ab$, to find a witness
$g\in\ab$ realizing $\alpha(f,g)$?  The Uniform Low Basis Theorem, proven
recently in \cite{BdBP12} by Brattka, de Brecht, and Pauly, will enhance this discussion dramatically.

The subgroup $\abz$ is countable and dense in $\ab$, as is every
subgroup $\abd$.  So likewise are the unions $\cup_{\bfd\in I}\abd$
over countable Turing ideals.  Consequently none of these groups
is closed in $\ab$, and none fixes any elements outside of $\Q$,
so the usual Galois correspondence between subfields of $\Qbar$
and closed subgroups of $\ab$ will not play a prominent role in this article.
Effective aspects of this correspondence were investigated by LaRoche
in \cite{L81}, and it seems likely that his results should relativize
to arbitrary Turing degrees.

It is fair to say that the results in this article are circumscribed by the boundaries
of pure computability theory:  to understand the proofs here, no significant
number-theoretic background is required.  Useful definitions and results
from computability are presented in the Appendix (Section \ref{sec:appendix}),
and \cite{S87} is the classic source for more basic background.  In ongoing work \cite{KM24},
the number theorist Debanjana Kundu has joined with the author to address
questions in this area where number theory plays an essential role,
and that role will likely only increase during any further investigations.
Certain relevant results from \cite{KM24} will be stated here without proof.

Additionally, the author would be remiss not to mention ongoing work by Wesley Calvert,
Valentina Harizanov, and Alexandra Shlapentokh.  While the goal of that work
is to develop a workable definition of a ``random'' algebraic field, the techniques
there use essentially the same presentation of $\Gal{\Qbar}{\Q}$ as given in
in Section \ref{sec:presentation}, and that work helped inspire the author
to consider $\Gal{\Qbar}{\Q}$ as a type-2-computable structure.
We hope that their results will appear soon.

\section{Presenting $\ab$}
\label{sec:presentation}

The celebrated theorem of Michael Rabin in \cite{R60} shows that, for every computable
field $F$, there is a computable presentation of the algebraic closure $\overline F$.
We immediately fix one computable presentation $\Qbar$ of the algebraic closure of $\Q$.
The central claim in Rabin's Theorem shows that $\Q$ itself forms a decidable subset of $\Qbar$,
since $\Q$ has a splitting algorithm, described by Kronecker in \cite{K1882}.  We also remark
that our choice of the specific presentation $\Qbar$ makes no difference:  the algebraic closure
of $\Q$ is computably categorical, as first shown by Ershov in \cite{E77}, meaning
that for any two computable presentations, there is a computable isomorphism between them,
allowing all results about either presentation to be transferred directly to the other.
For readers desiring more details, \cite{M21} is a useful beginning source.

Kronecker's theorems show that not only is $\Q$ decidable within $\Qbar$, but so is
every finite extension $F$ of $\Q$ within $\Qbar$, and the decision procedure for each such $F$
is uniform in each finite generating set for $F$.  Moreover, given a finite generating set,
we may effectively and uniformly find a primitive generator for $F$ (as seen in \cite{FJ86},
for instance).  Therefore, the following procedure is entirely computable
and defines an increasing nested sequence $F_0\subset F_1\subset\cdots$ of number fields,
all normal over $\Q$ (hence with each $F_{n+1}$ normal over $F_n$), such that $\cup_n F_n = \Qbar$.
Start with $F_0=\Q$, letting $z_0$ to be its multiplicative identity element.  Then,
for each $n\geq 0$, let $x_{n+1}$ be the least element of $\Qbar\setminus F_n$
(where ``least'' is defined using the usual order $<$ on the domain $\omega$ of $\Qbar$),
find all $\Q$-conjugates of $x_{n+1}$ in $\Qbar$, find a primitive generator $z_{n+1}$
of the field generated by all those conjugates (including $x_{n+1}$ itself), and define
$F_{n+1}=\Q(z_{n+1})$.

The absolute Galois group $\Gal{\Qbar}{\Q}$ is simply the automorphism group
$\ab$ of our field $\Qbar$, since every automorphism of $\Qbar$ fixes
$\Q$ pointwise.  We will refer to this group as $\ab$ from here on.
The \emph{Galois tree} $\TQbar$ will yield our representation of the elements of $\ab$.
Using the primitive generator $z_n$ of each field $F_n$, we define the elements
of the $n$-th level of $\TQbar$ to be the $n$-tuples
$$ (\sigma(z_1),\ldots,\sigma(z_n)) \in F_n^n,$$
where $\sigma$ ranges over $\an$.  This is computable:  the number of such $n$-tuples
is exactly the degree of $z_n$ over $\Q$, which is computable,
and their final coordinates are precisely the $\Q$-conjugates of $z_n$.
From the final coordinate $z$ it is easy to compute all previous coordinates:
$z_n$ generates $F_n$, so the map $z_n\mapsto z$ generates an automorphism
$\sigma$ of all of $F_n$, yielding the tuple.  Each initial segment of length $l<n$
has already appeared as a node at level $l$ in the tree. (Thus this is indeed a tree.)
We often refer to the individual node thus described as an automorphism
of $F_n$, thinking of it as the map $\sigma$ rather than
as the element $\sigma(z_n)$.  If $\tau\in\am$ extends $\sigma\in\an$, we will write
$\sigma\sqsubseteq\tau$.  Similarly, $\sigma\sqsubset f$ means that $f\in\ab$
extends $\sigma$.

The automorphisms of $\Qbar$ itself correspond bijectively to the paths through $\TQbar$,
in the obvious way:  each $f\in\ab$ may be viewed 
as the countable sequence $(f(z_1),f(z_2),\ldots)$, all of whose finite initial segments
are nodes in $\TQbar$.  Conversely, each path through $\TQbar$ names compatible
automorphisms of all fields $\F_n$, so their union is an automorphism of $\Qbar$.
The tree is a specific computable way of presenting these automorphisms,
by defining them on increasingly large number fields within $\Qbar$.  Moreover, the operations
of composition and inversion on $\ab$ are effective, in the sense that there exist
Turing functionals $\Gamma$ and $\Upsilon$ such that, for every $f,g\in\ab$ and every $n$,
$$ \Gamma^{(f(z_1),f(z_2),\ldots)\oplus(g(z_1),g(z_2),\ldots)}(z_n) = f(g(z_n))~~~\text{and}~~~
\Upsilon^{(f(z_1),f(z_2),\ldots)}(z_n) = f^{-1}(z_n).$$
Thus $\Gamma$ and $\Upsilon$ compute $f\circ g$ and $f^{-1}$ effectively from their oracles
$f\oplus g$ and $f$, respectively.  Indeed, each individual value $f(g(z_n))$ can be computed
just from $f(z_n)$ and $g(z_n)$, and $f^{-1}(z_n)$ can be computed just from $f(z_n)$.
(The \emph{join} $f\oplus g$ is simply the least upper bound
of $f$ and $g$ under Turing reducibility, often presented by splicing the two sequences
together as $(f(z_1),g(z_1),f(z_2),g(z_2),\ldots)$.  This join is not itself a path in the tree;
it is simply a convenient way to describe two paths simultaneously.)
All of this follows from work of Kronecker in \cite{K1882}; more modern
descriptions appear in \cite{E84} and \cite{M21}.

At times we will wish to consider several elements of $\ab$ simultaneously.
For this purpose we define $\TQbar^m$.  The level $n$ of $\TQbar^m$
consists of all $m$-tuples of nodes at level $n$ in $\TQbar$, with
$$ (\tau_1,\ldots,\tau_m) \sqsubseteq (\sigma_1,\ldots,\sigma_m) \iff
(\forall i\leq m)\tau_i\sqsubseteq\sigma_i.$$
\comment{
This is entirely natural, except for our treatment of the degenerate case $m=0$.
Every node in $\TQbar^0$ is the empty tuple, by definition, but we distinguish
these empty tuples according to the fields $F_n$ for which they are supposed
to list automorphisms.  That is, $\TQbar^0$ consists of a single node at each level $n$.
Its unique path, which contains every node of $\TQbar^0$, is the empty tuple $(\ab)^0$.
This convention will be convenient for us in certain situations.
}

For future use, we give a more general definition that covers both this situation
and that of other common presentations of uncountable structures.
\begin{defn}
\label{defn:treepres}
Let $\A$ be a structure of power at most the continuum.
A \emph{computable tree quotient presentation} of $\A$ consists of:
\begin{itemize}
\item
a computable subtree $T$ of $\omega^{<\omega}$; and
\item
a $\Pi^0_1$ equivalence relation $\sim$ on $[T]^2$
(where $[T]$ is the set of all paths through $T$); and
\item
for each $n$-ary function symbol $f$, a Turing functional $\Phi_f$ such that,
 whenever all $P_i$ lie in $[T]$, also $\Phi_f^{P_1\oplus\cdots\oplus P_n}\in[T]$; and
\item
for each $n$-ary relation symbol $R$, an effectively open set
$\U_R\subseteq [T]^n$; and
\item
for each constant symbol $c$, a computable path $P_c$ in $[T]$,
\end{itemize}
where all $\Phi_f$ and $\U_R$ respect $\sim$ and
$\A$ is isomorphic to the structure with domain $[T]/\!\sim$
under the functions $\Phi_f$, the relations $\U_R$, and the constants $[P_c]_{\sim}$.

Such an object is a \emph{computable tree presentation}
if $\sim$ is simply equality on $[T]$.
\end{defn}
Definition \ref{defn:treepres} naturally relativizes to the degree of any set $D$,
by allowing a $D$-oracle in each item.  If one wants a relation $R$
to be fully computable, this can be accomplished by putting both $R$ and
$\neg R$ in the signature.  (However, in case $T$ is finite-branching,
this leads to very simple relations!)  The usual presentation of the ordered field $\R$
is a computable tree quotient presentation:
$[T]$ is the class of all fast-converging Cauchy sequences of rational numbers, and
$\sim$ is the relation of having the same limit.
Notice that here $\sim$ is strictly $\Pi^0_1$, while $<$ is strictly $\Sigma^0_1$.
Automorphism groups of computable structures have computable tree presentations,
although often the tree must have terminal nodes.
Finally, countable structures $\A$ that satisfy the usual definition of a computable structure
(having domain $\omega$ and decidable atomic diagram) have natural
computable tree presentations:
just let $T$ contain every node of the form $(n,n,\ldots,n)$ with $n\in\omega$,
using the path $(n,n,n,\ldots)\in[T]$ to represent the element $n$ of the domain of $\A$.

Part of the reason for stating Definition \ref{defn:treepres} here is to emphasize what a simple structure
$\ab$ really is.  Our tree $\TQbar$ is finite-branching and has no terminal nodes.
(The finite-branching property makes it topologically compact, which
will be exploited dramatically in the rest of this article.)
There are no relation symbols, and the relation $\sim$ is just equality on $[\TQbar]$,
eliminating much of the uncertainty present in computable tree presentations
of other structures.

\section{Relating $\ab$ to $\an$}
\label{sec:Konig}

Having noted the simplicity of our presentation of $\ab$, we add that there remains some
undecidability about it:  the equality relation on its elements (i.e., on paths through $\TQbar$)
is only a $\Pi^0_1$ property, not decidable.  This situation informs our decision
to begin by examining positive formulas.  These will behave much more nicely than
non-positive formulas, as the avoidance of negation means that the formulas
will define $\Pi^0_1$-classes of automorphisms in $\ab$.

\begin{defn}
\label{defn:positive}
A formula $\phi$ in prenex form in a first-order language is \emph{positive} if $\phi$
does not use the negation connective $\neg$.
\end{defn}

Thus, in group theory, a positive sentence can say that there exists a group
element satisfying a particular disjunction of conjunctions of equations in words in other group elements,
but it cannot say that there exist two distinct such elements, as this would require negation.
(Recall that a \emph{sentence} is just a formula with no free variables.)

\begin{prop}
\label{prop:inf}
Let $\phi(\avec)$ be a positive sentence in the language of groups, with parameters
$\avec \in (\ab)^k$.  
Then
$$ \ab\models\phi(\avec)~~~\iff~~~(\forall n)~\an\models \phi(a_1\res F_n,\ldots,a_k\res F_n).$$
\end{prop}
Below we abbreviate $\phi(a_1\res F_n,\ldots,a_k\res F_n)$ by $\phi(\avec\res F_n)$.
\begin{pf}
We work by induction on the number of quantifiers.  A quantifier-free positive
sentence $\phi(\avec)$ may be taken to be a disjunction of conjunctions
of equations in $\avec$.  This is the base case for the induction.
If the disjunction holds in every $\an$,
then some single disjunct holds in infinitely many $\an$,
and it is then immediate that this disjunct holds in $\ab$.
Conversely, if the disjunction holds in $\ab$, then $\phi(\avec\res F_n)$ holds in $\an$
for every $n$, because when an equation holds in $\aut{F_{n+1}}$ (or in $\ab$),
it will remain true in $\an$ provided that we restrict the terms to $F_n$.
Here we see why the proposition requires $\phi$ to be positive:
inequations can be true in $\aut{F_{n+1}}$ yet false when restricted to $\an$.
(For instance, $\exists G_1\exists G_2~G_1\neq G_2$ holds in $\aut{F_1}$ but not in $\aut{F_0}$.)

For the inductive step, suppose first that $\phi(\avec)$ is of the form
$(\forall G)~\psi(G,\avec)$ with $\psi$ positive.  
Assume that $\ab\models (\forall G)~\psi(G,\avec)$.  Fix any $n$
and any $\sigma\in\an$, and choose an extension $g\sqsupset\sigma$
with $g\in\ab$.  Then $\ab\models \psi(g,\avec)$, and the inductive hypothesis
shows that $\an\models \psi(g\res F_n,\avec\res F_n)$.  Since $g\res F_n=\sigma$,
and since this works for every $n$ and $\sigma$, we have
$\an\models(\forall G)~\psi(G,\avec\res F_n)$ for all $n$.

Conversely, suppose that $\an\models(\forall G)~\psi(G,\avec\res F_n)$ for every $n$,
and pick any $g\in\ab$.  Then, for every $n$, $\an\models \psi(g\res F_n,\avec\res F_n)$.
The inductive hypothesis now shows that $\ab\models\psi(g,\avec)$, and thus
$\ab\models (\forall G)~\psi(G,\avec)$. as required.

Now consider a sentence $(\exists G)~\psi(G,\avec)$, with $\psi$ just as above.
We do the easy direction first,
assuming that this sentence holds in $\ab$, so that some $g\in\ab$ has
$\ab\models \psi(g,\avec)$.  By inductive hypothesis, we see that
for every $n$, $\an\models \psi(g\res F_n,\avec\res F_n)$,
hence $\an\models (\exists G)~\psi(G,\avec\res F_n)$ as required.

For the opposite direction, suppose that 
$\an\models (\exists G)~\psi(G,\avec\res F_n)$ for every $n$.
Now we use the tree $\TQbar$ deveoped earlier, building the subtree
$$ T_{\psi,\avec} = \set{\sigma\in\an}{n\in\omega~\&~\an\models\psi(\sigma,\avec\res F_n)}.$$
(Recall that the nodes of $\TQbar$ at level $n$ correspond to the automorphisms of $F_n$.)
Notice first that $T_{\psi,\avec}$ really is a subtree:  again, this holds because we require $\phi$
(and hence $\psi$) to be positive.  Negative sentences true in $F_{n+1}$ can fail in
$F_n$ when we consider their restrictions.

By assumption this $T_{\psi,\avec}$ has infinitely many nodes, and as a subtree of $\TQbar$,
it must be finite-branching.  K\"onig's Lemma therefore implies that it contains
an infinite path $\sigma_0\sqsubseteq\sigma_1\sqsubseteq\cdots$,
corresponding to an automorphism $g=\cup_n\sigma_n$ in $\ab$.
But since $\an\models\psi(\sigma_n,\avec\res F_n)$ for every $n$,
our inductive hypothesis shows that $\ab\models\psi(g,\avec)$ and therefore that
$\ab\models(\exists G)~\psi(G,\avec)$.
\qed\end{pf}

The final paragraph here begins to reveal the importance of K\"onig's Lemma,
which states that every infinite finite-branching tree has an infinite path.
Since we have presentations of all of the finite groups $\an$, uniformly in $n$,
the subtree $T_{\psi,\avec}$ of $\TQbar$ is decidable relative to the parameters $\avec$,
and the decision procedure is uniform in these parameters.  In Section \ref{sec:elementarity}
we will apply near-effective versions of K\"onig's Lemma in a similar way to derive further results.

Proposition \ref{prop:inf} has an obvious corollary, which surely was already folklore.

\begin{cor}
\label{cor:theory}
The positive theory of the group $\ab$ forms a $\Pi^0_1$ set.
\qed\end{cor}

\section{Scott ideals and partial elementarity}
\label{sec:elementarity}

Our ultimate goal is to find countable subgroups
of $\ab$, defined using computability theory, that are elementary subgroups
in the usual model-theoretic sense.  The question of full elementarity goes
beyond the scope of this article, but we will establish elementarity
at least for formulas of moderate complexity, provided that the
subgroup is defined by a Scott ideal.  Section \ref{subsec:lowbasis}
in the appendix recalls this notion.
The use of K\"onig's Lemma in the proof of Proposition \ref{prop:inf}
makes Scott ideals relevant to our topic.  

\begin{defn}
\label{defn:G_I}
For each Turing ideal in the Turing degrees, write
$$ \abI = \set{f\in\ab}{\deg{f}\in I}.$$
\end{defn}
It is immediate that $\abI$ is a subgroup of $\ab$, as arbitrary automorphisms $f$ and $g$
will have $f^{-1}$ and $f\circ g$ computable from $f$ and from the join $f\oplus g$ (respectively).
The closure of Turing ideals under join and downwards under $\leq$ does the rest.

The finite-branching subtree $T_{\psi,\avec}$ of $\TQbar$ that we built in
Proposition \ref{prop:inf} was computable from the tuple $\avec$, and its branching is
also $\avec$-computable.  Every path through it is a path through $\TQbar$
(hence an automorphism of $\Qbar$) realizing $\psi(G,\avec)$.
Conversely, every $g\in\ab$ with $\ab\models\psi(g,\avec)$ appears
as a path through $T_{\psi,\avec}$.  (Consequently, if $\ab\models\neg(\exists G)~\psi(G,\avec)$,
then $T_{\psi,\avec}$ is a finite tree.)  This will allow us to prove
the principal results of this section, Proposition \ref{prop:poselem} and
Theorem \ref{thm:compactlypos}.  As a warm-up, we begin with existential sentences,
which we no longer require to be positive.

\begin{thm}
\label{thm:Sigma1}
If $I$ is a Scott ideal in the Turing degrees,
then $\abI$ is \emph{elementary for purely existential sentences}.
By this we mean that whenever $\phi(\avec)$ is a sentence in the language of groups,
of the form $\exists G_1\cdots\exists G_m\psi(\Gvec,\avec)$
with $\psi$ quantifier-free and with parameters $\avec$ all in $\abI$,
$$  \ab\models \phi(\avec) \iff \abI \models\phi(\avec).$$
The same holds for purely universal sentences $\phi(\avec)$.
\end{thm}
\begin{pf}
Let $\phi(\avec)$ be $(\exists\Gvec)~\psi(\Gvec,\avec)$, with $\psi$ quantifier-free
and $\avec\in (\abI)^k$.  Again it is clear that a witness for $\phi$ in $\abI$ would
also be a witness in $\ab$.  For the difficult direction, suppose $\gvec\in(\ab)^m$
with $\ab\models\psi(\gvec,\avec)$. Viewing $\psi$ as $\vee\psi_i$, a disjunction
of sentences $\psi_i(\Gvec,\avec)$ each of which is a conjunction of equations
and inequations, assume without loss of generality that $\ab\models\psi_1(\gvec,\avec)$.
We divide $\psi_1$ into two (possibly empty) conjuncts $\psi_1^=(\gvec,\avec)$,
a conjunction of equations, and $\psi_1^{\neq}$, a conjunction of inequations.

By Proposition \ref{prop:inf}, every $n$ has $\an\models\psi_1^=(\gvec\res F_n,\avec\res F_n)$.
For each inequation in $\psi_1^{\neq}$, there exists some level $n$ such that
$\gvec\res F_n$ and $\avec\res F_n$ satisfy the inequation in $\an$.  
Moreover, satisfaction of inequations is preserved upwards in the tree.
Therefore, taking the maximum over all inequations in $\psi_1^{\neq}$
gives a level $n_0$ such that, for every $n\geq n_0$
and every $\vec{\sigma}\in (\an)^m$ with all $\sigma_i\res F_{n_0}=g_i\res F_{n_0}$
$\an\models\psi_1^{\neq}(\vec{\sigma},\avec)$.  In short,
above the node $\gvec\res F_{n_0}$ in the tree $\TQbar^m$,
$\psi_1^{\neq}$ is satisfied everywhere.  (We set $n_0=0$
in case $\psi_1^{\neq}$ is the empty conjunction.)

Fixing this $n_0$, we turn to the equations in $\psi_1^=$.  The procedure here
now recalls that of Proposition \ref{prop:poselem}, only beginning at the node $\gvec\res F_{n_0}$.
We define the set
$$ T_{\psi_1^=,\avec} = \set{\gammavec\in(\an)^m}{n\in\omega~\&~\gvec\res F_{n_0}\not\perp\gammavec
~\&~\an\models\psi_1^=(\gammavec,\avec\res F_n)}.$$
The second clause says that $\gammavec$ is compatible with $\gvec\res F_{n_0}$, so 
$T_{\psi_1^=,\avec} $ contains a single node $\gvec\res F_n$ at each level $n\leq n_0$
and then is allowed to start branching out after that.  
$T_{\psi_1^=,\avec}$ is a decidable subset of $(\TQbar)^m$
relative to an oracle for the tuple $\avec$ (along with its root node $\gvec\res F_{n_0}$,
which is finitely much information).  Proposition \ref{prop:inf} shows that for every $n\geq n_0$,
$\an\models\psi_1^=(\gvec\res F_n,\avec\res F_n)$, so $T_{\psi_1^=,\avec}$ must be infinite,
and its branching is $\avec$-computable (and finite!) just as before.  The Scott ideal $I$
must contain a degree $\bfd$ that is \PA\ relative to $\avec$, and this $\bfd$ must compute
a path $\gvec_I$ through $T_{\psi_1^=,\avec}$.
We may view $\gvec_I$ as an element of $(\ab)^m$, hence in $(\abI)^m$.  Since
every $n$ has $\an\models\psi_1^=(\gvec_I\res F_n,\avec\res F_n)$, we know that
$\ab\models\psi_1^=(\gvec_I,\avec)$, and since $\psi_1^=$ is quantifier-free,
$\psi_1^=(\gvec_I,\avec)$ holds in the subgroup $G_i$ as well.  
We also know that $\psi_1^{\neq}(\gvec_I,\avec)$ holds in $\abI$, because it holds at every node
above $\gvec\res F_{n_0}$.  Thus $G_i\models\psi_1(\gvec_I,\avec)$, and so $\abI\models\psi(\gvec_I,\avec)$
as desired.
\qed\end{pf}

\comment{
\begin{defn}
\label{defn:cpositive}
A formula $\psi$ in the language of groups (possibly with parameters $\avec$
from a particular group) is \emph{positive for compact spaces}
if $\psi$ is logically equivalent to a formula of the form
$$ (\exists G_1)\cdots(\exists G_m)(\forall H_1\cdots\forall H_l)~
[\alpha(\Gvec,\Hvec,\avec)~\&~\phi(\Gvec,\Hvec,\avec)],$$
where $\phi$ is positive (but allows quantifiers) and $\alpha$
is a conjunction of disjunctions of inequations.
\end{defn}
}

Now we wish to move our study of elementarity beyond existential and
universal sentences.  This becomes more difficult.  Our best result will be
Theorem \ref{thm:compactlypos}, which covers the class defined in Definition
\ref{defn:separated} below, containing more complex sentences,
including all positive sentences, but by no means all sentences.

\begin{defn}
\label{defn:QFneg}
A formula is \emph{quantifier-free negative} if it is
built from negated atomic formulas using $\wedge$ and $\vee$.
\end{defn}
Thus,  in the language of groups, every quantifier-free negative formula
is logically equivalent to a conjunction of disjunctions of inequations.

\begin{defn}
\label{defn:separated}
A formula $\beta$ is \emph{separated} (and is also
\emph{$\Sigma_0$-separated} and \emph{$\Pi_0$-separated})
if either of the following holds.
\begin{itemize}
\item
$\beta$ itself is positive; or
\item
$\beta$ is the conjunction of a positive formula and a quantifier-free negative formula.
(The positive formula is allowed to use quantifiers.)
\end{itemize}
A formula $\alpha$ is \emph{$\Sigma_{n+1}$-separated} if $\alpha$
is of the form $\exists\Gvec\beta(\Gvec)$, where $\beta$ is
$\Pi_n$-separated; while $\alpha$ is \emph{$\Pi_{n+1}$-separated} if it is of the form
$\forall\Gvec\beta(\Gvec)$, where $\beta$ is $\Sigma_n$-separated.
\end{defn}
Notice that negations of separated formulas need not be separated.
More generally, the negation of a $\Sigma_n$-separated formula
need not be $\Pi_n$-separated, nor vice versa.
Also, a quantifier-free formula in a form such as $(A\vee(\neg B))\wedge(C\vee(\neg D))$
(where $A$, $B$, $C$, and $D$ are equations) will be neither.
It is routine to show that every positive formula is logically equivalent to a positive formula
in prenex normal form.  However, negation does not move readily across quantifiers,
and so we make is no similar claim for separated formulas.

Now we begin with Proposition \ref{prop:poselem}, concerning positive sentences,
which we will then extend to Theorem \ref{thm:compactlypos}
for $\Sigma_2$-separated sentences.

\begin{prop}
\label{prop:poselem}
Let $I$ be any Scott ideal in the Turing degrees.
Then $\abI$ is an \emph{elementary subgroup of $\ab$ for positive sentences}.
By this we mean that whenever $\phi(\avec)$ is a positive sentence
in the language of groups, with parameters $\avec$ all in $\abI$,
$$  \ab\models \phi(\avec) \iff \abI \models\phi(\avec).$$
(The same therefore holds when $\phi(\avec)$ is the negation of a positive sentence,)
\end{prop}
This is just the standard notion of an elementary substructure, restricted
to positive sentences. The question of whether this ``positive elementarity'' holds
for all Turing ideals (not just Scott ideals) will appear among our open questions
in Section \ref{sec:questions}.
\begin{pf}
As usual for questions of elementarity, the important step is to
``go downwards'' for $\exists$-sentences.  The equivalence in the theorem
is immediate in case $\phi$ is quantifier-free.  Moreover, once it is known to hold for
sentences with $m$ quantifiers, the case of a sentence $(\exists G)~\psi(G,\avec)$
(where $\psi$ has $m$ quantifiers) that holds in $\abI$ is quickly handled, as the witness
$g\in \abI$ with $\abI\models\psi(g,\avec)$ also satisfies $\ab\models\psi(g,\avec)$
by induction.  The converse is the step requiring work.

Assume therefore that $\ab\models(\exists G)~\psi(G,\avec)$, and let $g\in\ab$
witness this fact.  Of course, $g$ need not lie in $\abI$.  As in Proposition \ref{prop:inf},
we define the subtree
$$ T_{\psi,\avec} = \set{\sigma\in\an}{n\in\omega~\&~\an\models\psi(\sigma,\avec\res F_n)}.$$
By Proposition \ref{prop:inf}, for every $n$, $g\res F_n$ lies in $T_{\psi,\avec}$.
(Thus, as we remarked above, every witness is given by a path in $\TQbar$ which
lies entirely within $T_{\psi,\avec}$.)  It follows that $T_{\psi,\avec}$ is infinite.
Moreover, since every $\an$ is finite and effectively presented (uniformly in $n$),
determining whether a given $\sigma\in\TQbar$ lies in $T_{\psi,\avec}$
only requires us to know the elements $a_1\res F_n,\ldots,a_k\res F_n$ within $\an$.
This merely requires an oracle (formally, $a_1\oplus\cdots\oplus a_k$) for $\avec$,
so $T_{\psi,\avec}$ is an $\avec$-decidable subtree of $\TQbar$.  Moreover,
as the branching in $\TQbar$ is computable, with $\avec$ we can
now compute the number of immediate successors of any
$\sigma\in T_{\psi,\avec}$:  just compute the number of immediate successors of
$\sigma$ in $\TQbar$, find them all, and use $\avec$ to check which ones lie in $T_{\psi,\avec}$.
(The answer could be zero, of course:  $T_{\psi,\avec}$ may have terminal nodes.)

As $I$ is a Turing ideal, the join $\deg{a_1}\cup\cdots\cup\deg{a_k}$ lies in $I$.
Being a Scott ideal, $I$ therefore contains some degree $\bfd$ that is \PA\ relative
to $\deg{a_1}\cup\cdots\cup\deg{a_k}$.  By Definition \ref{defn:PA}, $\bfd$
can compute a path $g_I$ through the infinite finite-branching subtree $T_{\psi,\avec}$.
This $g_I$ is therefore a path through $\TQbar$ as well, i.e., an automorphism of $\Qbar$,
and lies in $\abI$ because it is computable from $\bfd$, which lies in $I$.
By the definition of $T_{\psi,\avec}$, every $n$ has $\an\models\psi(g_I\res F_n,\avec\res F_n)$.
Proposition \ref{prop:inf} now shows that $\ab\models \psi(g_I,\avec)$, so
by inductive hypothesis $\abI\models\psi(g_I,\avec)$ as well,
so $\abI\models (\exists G)~\psi(G,\avec)$ as required.
\comment{
Finally we remark that if $\phi(\avec)$ is of the form $(\forall G)~[\neg\psi(G,\avec)]$
where $\psi$ is positive and has $m$ quantifiers, then the equivalence in the theorem
is now known to hold of its negation $(\exists G)\psi(G,\avec)$,
and this immediately implies the same equivalence for $\phi(\avec)$ itself.
This completes the induction on $m$.
}
\qed\end{pf}

\begin{thm}
\label{thm:compactlypos}
Let $I$ be any Scott ideal in the Turing degrees. Then the subgroup
$\abI$ of $\ab$ is elementary for $\Sigma_2$-separated sentences,
and also for their negations.
This means that every such sentence (with parameters from $\abI$)
holds in $\abI$ just if it holds in $\ab$.
\comment{
Moreover, for any $A\subseteq\omega$, let $I_A=\set{\bfd}{\bfd\leq\bfb~\&~\bfd\leq\bfc}$
be the specific Scott ideal defined in Theorem \ref{thm:lowexactpair}.
Then, uniformly in $A$, we can compute an approximation to the jump
of a witness to the truth of each compactly-positive sentence $\psi$.
}
\end{thm}
\begin{pf}
Definition \ref{defn:separated} gives the form of $\psi(\avec)$ in the important case:
$$ (\exists G_1\cdots\exists G_m)(\forall H_1\cdots\forall H_l)~
[\alpha(\Gvec,\Hvec,\avec)~\&~\phi(\Gvec,\Hvec,\avec)],$$
where $\alpha$ is quantifier-free negative and $\phi$ is positive.
(For the other case, where $\psi(\avec)$ itself is positive,
apply Proposition \ref{prop:poselem}.)

First suppose $\psi(\avec)$ holds in $\abI$,
with some specific witness tuple $\gvec_0$ there.
By Theorem \ref{thm:Sigma1}, the $\Pi_1$ sentence $\forall\Hvec\alpha(\gvec_0,\Hvec,\avec)$
will also hold in $\ab$, as will the positive sentence $\forall\Hvec\psi(\gvec_0,\Hvec,\avec)$,
by Proposition \ref{prop:poselem}.
Thus $\ab\models\psi(\avec)$, as required by the theorem.

\comment{
we handle the easy direction.  If $\abI\models \psi(\avec)$,
then fix a tuple $\gvec_I$ of witnesses in $\abI$, so
$\abI\models \forall\Hvec[\alpha(\gvec_I,\Hvec,\avec)~\&~\phi(\gvec_I,\Hvec,\avec)]$.
If there were an $\hvec$ in $\ab$ such that
$\ab\models\neg [\alpha(\gvec_I,\hvec,\avec)~\&~\phi(\gvec_I,\hvec,\avec)]$,
then for some $n\in\omega$ we would have
$\an\models \neg \alpha(\gvec_I\res F_n,\hvec\res F_n,\avec\res F_n)$.
But then
$$ \set{\etavec\in(\TQbar)^l}{\etavec\sqsubseteq\hvec\text{~or~}
[\hvec\sqsubseteq\etavec~\&~\aut{F_{|\etavec|}}\models\phi(\gvec_I\res F_{|\etavec|},
\etavec,\avec\res F_{|\etavec|})]}$$
would be an infinite $(\gvec_I\oplus\avec)$-computable subtree of $(\TQbar)^l$,
so the Scott ideal $I$ would contain the degree of an infinite path $\hvec_I$ in this tree,
for which we would have $\abI\models \neg [\alpha(\gvec_I,\hvec_I,\avec)~\&~\phi(\gvec_I,\hvec_I,\avec)]$,
contradicting the assumption that 
$\abI\models \forall\Hvec[\alpha(\gvec_I,\Hvec,\avec)~\&~\phi(\gvec_I,\Hvec,\avec)]$.
Therefore no such $\hvec$ exists in $\ab$, and thus $\gvec_I$ witnesses
that $\ab\models\psi(\avec)$ as required.
}

For the converse, suppose that $\psi(\avec)$ holds in $\ab$.  
We need to establish that $\psi(\avec)$ holds in $\abI$ as well.
Fix a tuple $\gvec\in (\ab)^m$ of witnesses to its truth in $\ab$.
Defining $A=\oplus_i a_i$ and $\bfa=\deg{A}\in I$, we will produce an $\bfa$-computable
subtree $T_{\psi,\avec}$ of $(\TQbar)^m$, such that every path through
this tree is a tuple of witnesses to the truth in $\abI$ of (this equivalent of) $\psi(\avec)$.
Since the Scott ideal $I$ must contain the degree of some such path, this will
prove the Corollary.

The twist here is that we give an $\bfa$-computable approximation of
the subtree $T_{\psi,\avec}$ of $\TQbar$, rather than a decision procedure.
There may be finitely many stages $s+1$ at which
$T_{\psi,\avec,s+1}$ deletes nodes from $T_{\psi,\avec,s}$.
However, from some stage $s_0$ on, all approximations will define the same subtree of $\TQbar$,
and so this subtree $T_{\psi,\avec}$ will indeed be an $\bfa$-computable subtree
of $(\TQbar)^m$.  This twist means that our construction of $T_{\psi,\avec}$
is not uniform in $\psi$ and $\avec$, as uniformity would also require
knowledge of the stage $s_0$.  Presenting the contruction this way
will facilitate a proof in Section \ref{sec:Skolem}.

At each stage $s$, we define $T_{\psi,\avec,s}$ as follows.
Consider the shortest $\gamma\in\TQbar$
such that both
$$ (\exists\delta\in(\TQbar)^m)~[|\delta|=s~\&~\gamma\not\perp\delta~\&~
\as\models \forall\Hvec~\phi(\delta,\Hvec,\avec)]$$
(which is $\bfa$-decidable, as level $s$ of $(\TQbar)^m$ is finite)
and
$$ \aut{F_{|\gamma |}} \models (\forall\Hvec)~\alpha(\gamma,\Hvec,\avec).$$
We will see below that such a $\gamma$ must exist.
Define $\gamma_s$ to be the leftmost such $\gamma$
(at this least possible level).
Define the $n$-th level of the tree $T_{\psi,\avec,s}$ to contain:
\begin{itemize}
\item
only $\gamma_s\res n$, if $n\leq |\gamma_s|$;
\item
those $\delta\in(\TQbar)^m$ of length $n$ such that
$$\gamma_s\sqsubseteq\delta~\&~
\an\models (\forall H_1\cdots\forall H_l)~\phi(\delta,\Hvec,\avec),$$
if $n>|\gamma_s|$ and this set is nonempty; or
\item
all $\delta\in(\TQbar)^m$ of length $n$ such that
$\delta\res(n-1)\in T_{\psi,\avec,s}$, if $n>|\gamma_s|$ and the set
in the preceding item is empty.
\end{itemize}
Thus $T_{\psi,\avec,s}$ contains $\gamma_s$ but no other node at level $s$.
Above that level, it contains all extensions of $\gamma_s$ that appear
to witness satisfaction of our formulas.  (However, it is possible that
the tree of all such extensions is finite; in this case, we recognize its finiteness
when we reach its maximal node(s), say of length $n_0$, and simply define it trivially from there on up,
so as to ensure that $T_{\psi,\avec,s}$ is infinite.  If this occurs,
then when we come to define $\gamma_{n_0+1}$, it will be different from $\gamma_s$.)
This completes stage $s$.

First we claim that the sequence $\la\gamma_s\ra_{s\in\omega}$ stabilizes
on a limit $\gamma$, which will be a specific node in $\TQbar$.
Indeed, for our witness $\gvec\in(\ab)^m$, we know first of all
that $\an\models (\forall\Hvec)\phi(\gvec\res F_n,\Hvec,\avec\res F_n)$ for every $n$,
by Proposition \ref{prop:inf}, and second of all that
$\ab\models(\forall\Hvec)\alpha(\gvec,\Hvec,\avec)$.
The quantifier-free formula $\alpha(\Gvec,\Hvec,\avec)$ is a finite conjunction
$\wedge_i\alpha_i$, where each $\alpha_i$ is a finite disjunction
of inequations.  Since $\ab\models(\forall\Hvec)\alpha_i(\gvec,\Hvec,\avec)$
and every inequation is realized in some $\an$,
$(\ab)^l$ is covered by the open sets
$$\set{\hvec\in (\TQbar)^l}{
\an\models\alpha_i(\gvec\res F_n,\hvec\res F_n,\avec\res F_n)}$$
as $n$ ranges over $\omega$.  The compactness of $(\ab)^l$ then
shows that there is some finite level $n_i$ with
$\aut{F_{n_i}}\models(\forall\Hvec)\alpha_i(\gvec\res F_{n_i},\Hvec,\avec\res F_{n_i})$.
Thus, for every $n\geq\max_i(n_i)$,
$$ \an\models (\forall\Hvec)\alpha(\gvec,\Hvec,\avec),$$
as the inequations in every $\alpha_i$ are all satisfied for every $\Hvec$
by that level.  Thus, at every stage $s\geq\max_i(n_i)$,
the restriction $\gvec\res F_{\max(n_i)}$ is a node which satisfies
the conditions for $\gamma_s$ (with $\gvec\res F_s$ as $\delta$).

This shows first that $\gamma_s$ is defined at every stage $s$,
as $\gvec\res F_{\max(n_i)}$ satisfies all the requisite conditions at every stage.
Second, whenever $\gamma_s\neq\gamma_{s+1}$,
$\gamma_s$ will never again be chosen as $\gamma_t$ for any $t>s$.
(The reason for it not to be chosen as $\gamma_{s+1}$ must stem from the lack
of an extension $\delta$ at level $s+1$, which will continue at subsequent stages $t$.)
Thus, among the finitely many nodes that are either shorter than 
$\gvec\res F_{\max_i(n_i)}$ or to its left at the level $\max_i(n_i)$,
each one is chosen as $\gamma_s$ either finitely often or cofinitely often.
It follows that $\la\gamma_s\ra_{s\in\omega}$ does indeed stabilize on a limit $\gamma$.

Let $s_0$ be a stage such that for every $s\geq s_0$,
$\gamma_s=\gamma$ is the limiting node.  Then for every $s\geq s_0$,
$T_{\psi,\avec,s}$ contains precisely those $\delta\sqsupseteq\gamma$ for which
$$ \aut{F_{|\delta|}} \models (\forall H_1\cdots\forall H_l)~
[\alpha(\delta,\Hvec,\avec\res F_{|\delta|})~\&~\phi(\delta,\Hvec,\avec\res F_{|\delta|})].$$
(Every level $n$ contains such a $\delta$,
since $\gvec\res F_n$ is one such.)  Therefore
$T_{\psi,\avec,s} = T_{\psi,\avec,s_0}$ for all $s\geq s_0$.
This shows that $\abI\models\psi(\avec)$, because
this $\bfa$-computable tree must contain a path $\gvec_0$
with Turing degree in the Scott ideal $I$, which will satisfy
$\aut{F_{|\gamma|}}\models \forall\Hvec\alpha(\gvec_0\res F_{|\gamma|},\Hvec,\avec\res F_{|\gamma|})$
and (for all $n$) $\an\models (\forall H_1\cdots\forall H_l)~\phi(\gvec_0\res F_n,\Hvec,\avec\res F_{|\gamma|})$.

The result for negations of $\Sigma_2$-separated sentences follows immediately
from that for $\Sigma_2$-separated ones.
\qed\end{pf}

The reader may wonder why we bothered with the initial stages
of this construction, rather than just choosing $\gamma=\gamma_{s_0}$
and beginning nonuniformly at stage $s_0$.  Doing so would indeed prove Theorem
\ref{thm:compactlypos}.  However, in the proof of Theorem
\ref{thm:Skolem} below, we will use the tree $T_{\psi,\avec}$
again and will want to know its construction uniformly through
all stages, even though we will still not know which stage is stage $s_0$.
It seems best to give the entire construction here, rather than revising
it during that proof.

\section{Definable sets for $\ab$}
\label{sec:definability}

Section \ref{sec:elementarity} concerned the extent to which specific
countable subgroups of $\ab$ can be elementary.  However, the methods
derived there also yield facts about $\ab$ itself.  In this section
we determine the complexity levels of certain definable subsets of $\ab$,
viewed as sets of reals within the hyperarithmetical hierarchy.
In Section \ref{sec:Skolem} the same methods will allow us to examine
the difficulty of computing Skolem functions in $\ab$ for various formulas.

In Proposition \ref{prop:poselem}, given a positive formula $\psi(\Gvec,\Avec)$
and a tuple $\avec$ from $\ab$, we produced an $\avec$-computable
subtree $T_{\psi,\avec}$ of $(\TQbar)^m$ such that the paths through
$T_{\psi,\avec}$ are precisely the $m$-tuples $\gvec\in (\ab)^m$ such that
$\ab\models \psi(\gvec,\avec)$.  In that Proposition, the point of the subtree
was that it must contain a path of PA-degree relative to $\avec$.  However,
the same construction also yields a description of the subset of $(\ab)^m$
defined by the formula $\psi(\Gvec,\avec)$ with the parameters $\avec$.
Thus we infer the following result.

\begin{thm}
\label{thm:posdefine}
For every positive formula $\psi(\Gvec,\Avec)$ in the language of groups,
and every tuple $\avec\in(\ab)^l$ of parameters, the subset $S_{\psi,\avec}$ of $(\ab)^m$ defined
by $\psi(\Gvec,\avec)$ is arithmetically a $\Pi^0_1$ set of reals relative to $\avec$,
and the $\Pi^0_1$ property defining it is uniform in both $\psi$ and $\avec$.
\end{thm}
\begin{pf}
The tuple $\gvec$ lies in $S_{\psi,\avec}$ just if, for every $n$, $\gvec\res F_n\in T_{\psi,\avec}$,
where this is the subtree of $(\TQbar)^m$ defined in the proof of Proposition \ref{prop:poselem}.
\qed/\end{pf}

Notice that, if the positive formula $\psi$ is $\Sigma_n$ in complexity,
then the definition it yields for the set $S_{\psi,\avec}$ is boldface
$\bf{\Sigma^1_n}$, as all quantification is over elements of $\ab$, which is to say,
over reals, not over natural numbers.  (A parameter-free definition $\psi$
would define a lightface $\Sigma^1_n$ set, but any parameters must be included
as an oracle.)  On its face, therefore, $S_{\psi,\avec}$ did not appear to belong to the
hyperarithmetical hierarchy at all unless $n=0$.  Theorem \ref{thm:posdefine}
shows $S_{\psi,\avec}$ to have far lower complexity than the definition would have suggested.

On the other hand, the ``warm-up case'' of Theorem \ref{thm:Sigma1} does not yield
this same result, nor does the construction in the proof of Theorem \ref{thm:compactlypos}.
These two theorems each allowed some use of negation in the defining formulas.
Now even the quantifier-free statement $G_0\neq G_1$ has arithmetic complexity $\Sigma^0_1$
when applied to $\ab$:  $g_0\neq g_1$ just if there exists some $n$ with $g_0\res F_n\neq g_1\res F_n$.
So, when the formulas in these theorems are not positive, there is no reason to expect
the sets they define to be $\Pi^0_1$-classes, even relative to the parameters.
In the proof of Theorem \ref{thm:Sigma1}, the tree we constructed was not uniform in the formula:
it required the knowledge of a level $n_0$ at which the inequations had all been satisfied.
Likewise, although we gave a uniform construction in the proof of Theorem \ref{thm:compactlypos},
that construction did not actually decide an $\avec$-computable subtree of $\TQbar$.
Instead, it approximated that subtree $T_{\psi,\avec}$, in such a way that the entire approximation
converged at some finite stage, but without giving any way to determine that finite stage effectively.
Indeed, in these cases, the subsets of $(\ab)^m$ defined by the formulas in question
(with parameters $\avec$) are arithmetically only $\Sigma^0_2$ relative to $\avec$.
(We leave the construction of the $\Sigma^0_2$ property to the reader.)
The $\Sigma^0_2$ property is uniform in both $\psi$ and $\avec$, once again,
and of course $\Sigma^0_2$ is still a substantial improvement over the original
(non-hyperarithmetical) definitions of these definable sets, but it is no longer clear that they are
$\Pi^0_1$-classes.  However, it remains open whether this is sharp, in the specific case of $\ab$:
does there actually exist an existential or universal formula, or a $\Sigma^0_2$-separated formula,
defining a subset of $(\ab)^m$ that is not a $\Pi^0_1$-class relative to the parameters?
We conjecture that the answer is affirmative, but that answer seems likely to require
a fair dose of number theory in the construction.

The uniformity in the constructions above yields an application to the theory of
$\ab$ as well.  If $\psi$ is a $\Sigma_1$ sentence in the language of groups,
say $(\exists\Gvec)\alpha(\Gvec)$ with $\alpha$ in disjunctive normal form,
then $\psi$ holds just if the decidable subtree $T_\psi$ of $\TQbar$ contains
a node $\sigma$ at which all inequations (in one or another disjunct of $\alpha$) hold
and such that every level $>|\sigma|$ contains an extension of $\sigma$ in which
all equations from that disjunct hold.  This is arithmetically a $\Sigma^0_2$ statement.
However, a $\Pi_1$ sentence $(\forall\Gvec)\beta(\Gvec)$, now with $\beta$
in conjunctive normal form, holds just if there is a level $l$ such that
for every $\sigma$ at level $l$ in $\TQbar$ there is a conjunct $\beta_i$ of $\beta$
such that every level $\geq l$ contains a successor of $\sigma$ satisfying the same
conjunct $\beta_i$.  Again this is arithmetically a $\Sigma^0_2$ property, so we have
extended Corollary \ref{cor:theory} to a stronger result.  We remind the reader
that the existential theory of $\ab$ is, on its face, $\Sigma^1_1$, as it quantifies over reals,
yet its complexity is far lower.
\begin{prop}
\label{prop:theory}
The existential theory of $\ab$ is $\emptyset'$-decidable, as it and its complement
are both $\Sigma^0_2$.
\qed\end{prop}

\section{Skolem functions for $\ab$}
\label{sec:Skolem}

\comment{
Sets definable by positive formulas without parameters are always closed,
and therefore are always describable as the set of paths through a decidable subtree
of $\ab$ (i.e., as a $\Pi^0_1$-class).  More generally, each set definable
by a positive formula with parameters $\fvec$ is the set of paths
through a subtree decidable relative to the join $\oplus\fvec$ of those parameters.
}

Consider a $\Pi_2$ formula $\forall F \exists G\alpha(F,G)$, with $\alpha$
quantifier-free, that holds in $\ab$.  A \emph{Skolem function} for this $\alpha$
is a function $\G:\ab\to\ab$ such that, for every $f\in\ab$, $\alpha(f,\G(f))$
holds in $\ab$.  Thus $\alpha$ is a kind of choice function:  knowing that
(for its input $f$) some $g$ satisfying $\alpha(f,g)$ must exist, it picks out
one such $g$.  By the Axiom of Choice, such a Skolem function must exist,
and in general there will be many Skolem functions for this given $\alpha$.
(Of course there might be only one:  this holds if the $g$ corresponding to each
$f$ is unique, so $\ab\models \forall F~\exists!G~\alpha(F,G)$.)
More generally, in any formula in prenex form, we can replace the $\exists$-quantifiers
by Skolem functions of all preceding variables, so that $\forall a\exists b\forall c\exists d R(a,b,c,d)$
becomes $\forall a\forall c~R(a,\F(a),c,\G(a,c))$, for example.

The meaning of ``Skolem function'' for the initial $\exists$-quantifier in a
$\Sigma_n$ sentence can be ambiguous:  what is required is essentially
a $0$-ary function, often viewed simply as a constant $a_0\in\ab$.  (That is,
$\exists a \forall b\exists c\forall d~R(a,b,c,d)$ Skolemizes to
$\forall b\forall d~R(a_0,b,\F(b),d)$.)  We will consider here
a more general situation:  even if $\ab\not\models \forall F\exists G\alpha(F,G)$,
we would like to have a function $\G$ such that every $f\in\ab$ for which
$\ab\models\exists G\alpha(f,G)$ has $\ab\models\alpha(f,\G(f))$.  In other words,
$\G$ should output a witness to $\exists G\alpha(f,G)$ whenever one exists.
We will refer to such a $\G$ as a \emph{generalized Skolem function}
for $\exists\Gvec\alpha(\Fvec,\Gvec)$, and ideally we would want to produce
such a $\G$ uniformly for each formula.

Notice that, if we have a way of producing these generalized Skolem functions
uniformly for $\Sigma_n$ formulas, then we can also produce
Skolem functions uniformly for $\Pi_{n+1}$ formulas.  Indeed,
if $(\forall\Fvec)(\exists\Gvec)\psi(\Fvec,\Gvec,\avec)$ is $\Pi_{n+1}$ with parameters $\avec$
(so $\psi$ is $\Pi_{n-1}$), then for every tuple $\fvec$, a generalized Skolem function $\G$ for
$\Sigma_n$ formulas can be applied to the formula $(\exists\Gvec)\psi(\fvec,\Gvec,\avec)$,
in which we consider $\fvec$ and $\avec$ together as one long tuple of parameters.
We see that, whenever this formula holds in $\ab$, we will have
$\ab\models\psi(\fvec,\G(\psi,\fvec,\avec),\avec)$, and so, using this $\G$,
we can compute $\G^*$ such that $\ab\models(\forall\Fvec)\psi(\Fvec,\G^*(\Fvec,\avec),\avec))$,
making $\G^*$ a Skolem function for the $\Pi_{n+1}$ formula
$(\forall\Fvec)(\exists\Gvec)\psi(\Fvec,\Gvec,\avec)$ with which we started.
Finally, finding this $\G^*$ was uniform in $\psi$ and $\avec$,
so this procedure succeeds uniformly for all $\Pi_{n+1}$ formulas.
Therefore, in theorems such as Proposition \ref{prop:Skolem}
and Theorem \ref{thm:Skolem} below, we will focus on $\Sigma_n$ formulas.

The existence of generalized Skolem functions again follows directly from the Axiom
of Choice.  However, as in the rest of this article, we are concerned with effectiveness
rather than mere existence.  For purposes of computation in $\ab$,
it would be ideal for Skolem functions to be computable,
and moreover they should even be uniform in the formula in question.
Thus, at the $\Sigma_1$ level, one would hope for a computable function
$\G(\alpha,\Fvec)$ of both the formula $\alpha(F_1,\ldots,F_{l_\alpha},\Gvec)$ and the parameters $\Fvec$,
satisfying for every quantifier-free $\alpha$ and every $\fvec\in(\ab)^{l_\alpha}$
$$ \ab\models(\exists G_1\cdots\exists G_m)\alpha(\fvec,\Gvec)~~~\implies~~~
\ab\models\alpha(\fvec, \G(\alpha,\fvec)).$$
If this held, we would then move on to more complex formulas and
pose similar questions about computing their Skolem functions.
If there is no such $\G$ -- that is,
if generalized Skolem functions for $\Sigma_1$ formulas in $\ab$
are not computable -- then before moving on to $\Sigma_2$ formulas,
we would next ask how close we can come to computing generalized
$\Sigma_1$ Skolem functions.

It turns out that the second of these situations is the case.
In the article \cite{KM24}, written in parallel to this article,
Kundu and the present author have established the following theorem.
\begin{thm}[Kundu \& Miller \cite{KM24}]
\label{thm:KM}
There is no computable generalized Skolem function for the formula
$(\exists G)G\circ G=F$.
Indeed, for every Turing functional $\Phi$, there exists an automorphism $f\in\ab$
such that $\ab\models(\exists G)~G\circ G=f$, yet either $(\Phi^f)\circ(\Phi^f)\neq f$
or else $\Phi^f\notin\ab$ (including the possibility that $\Phi^f$ is not even total).
In fact, such a counterexample $f$ may be computed uniformly from an index for $\Phi$.
\qed\end{thm}
Since it is impossible to compute a Skolem function for this specific formula,
it is certainly impossible to compute a generalized Skolem function for
$\Sigma_1$ formulas collectively.  Therefore we turn here to the second question:
how close can we come to computing Skolem functions?  It turns out that
we can come quite close.  To quantify this, consider the following concept.
\begin{defn}
\label{defn:superapproximable}
A function $F:2^\omega\to 2^\omega$ is \emph{superapproximable}
if there exists a Turing functional $\Phi$ such that, for every $A\in 2^\omega$,
$$ (\forall n\in\omega)~\lim_s \Phi^A(n,s) = \chi_{(F(A))'}(n).$$
If this holds, we say that $\Phi$ \emph{superapproximates} $F$.
The same definition applies to functions on $\ab$ (viewed as paths through $\TQbar$),
or on Baire space $\omega^\omega$.
\end{defn}
A superapproximation is generally not a computation of $F$,
but it comes quite close to being one.
It is stronger than an ordinary approximation, in which
the characteristic function $\chi_{F(A)}$ would be presented as the limit
of an $A$-computable function and hence would have $F(A)\leq_T A'$.
With a superapproximation, we approximate not just $F(A)$ but its
jump $(F(A))'$, uniformly in $A$, showing that $F(A)$ itself is low
relative to $A$.  From this it is not difficult to create
an approximation of $F(A)$ itself, if this is what one wants,
as there exists a single computable function $f$ that is a $1$-reduction
from $S$ to $S'$ simultaneously for every $S\subseteq\omega$.
(For the program for $f$, see \cite[Lemma 1]{M23}.)
However, the specific point is to approximate the jump $(F(A))'$
computably and uniformly in $A$, thus proving $F(A)$ to be uniformly low relative to $A$.

Subsection \ref{subsec:lowbasis} of the Appendix presents a version of the Uniform
Low Basis Theorem of Brattka, de Brecht, and Pauly as Theorem \ref{thm:ULB}.
It follows from their theorem that there is a superapproximable function $F$
such that, for every $A\in 2^\omega$, $F(A)$ has PA-degree relative to $A$.

\begin{lemma}
\label{lemma:super^2}
Suppose $F,G:2^\omega\to 2^\omega$ are superapproximable.
Then so is $G\circ F$, and a superapproximation for $G\circ F$ may be found
uniformly in those for $G$ and $F$.
\end{lemma}
\begin{pf}
Fix $A\in 2^\omega$, set $B=F(A)$ and $C=G(B)$, and
let $\Phi$ and $\Gamma$ be functionals such that, for every $e$,
$$ \lim_s \Phi^A(e,s)=\chi_{B'}(e)~~~\&~~~\lim_s \Gamma^B(e,s)=\chi_{C'}(e).$$
From $\Phi$ we derive, uniformly, a Turing reduction $B'\leq_T A'$.
(For each $e$, just use $A'$ to find a stage $s_e$ after which
$\Phi^A(e,s)$ is constant.)  Similarly from $\Gamma$ we get $C'\leq_T B'$ uniformly.
Composing Turing reductions is a uniform process, so we have a reduction
$C'\leq_T A'$ uniformly, say $\chi_{C'}=\Theta^{A'}$.  Also, uniformly in $A$,
we have an $A$-computable enumeration $\la A'_s\ra_{s\in\omega}$ of $A'$,
Now we define $ \Psi^A (e,s)$  to to output the value $\Theta_u^{A'_t}(e)$,
where $\la t,u\ra$ is the least pair such that $t\geq s$ and
$\Theta_u^{A'_t}(e)\converges$.  Such a pair must exist,
because $\Theta^{A'}(e)\converges$.  Moreover,
$\Theta^{A'_t}(e)=\Theta^{A'}(e)$ for all sufficiently large $t$,
so for all sufficiently large $s$ we have $\Psi^A(e,s)=\chi_{C'}(e)$.
\qed\end{pf}

With this notion, we now address the question of the complexity of Skolem
functions for various formulas about $\ab$.  In large part, the answers
will be derived by exploiting the uniformity of our earlier results on elementarity.

\begin{prop}
\label{prop:Skolem}
$\ab$ has superapproximable generalized Skolem functions for every
positive $\Sigma_n$ formula $\exists G\psi(\avec,G)$, for every $n$.
Moreover, superapproximations of these functions can be given uniformly
in $n$ and the formula $\psi$.
\end{prop}
As remarked above, it follows that $\ab$ has uniformly superapproximable
Skolem functions for all positive $\Pi_n$ formulas as well.
\begin{pf}
The same tree $T_{\psi,\avec}$ used in the proofs of Propositions \ref{prop:inf}
and \ref{prop:poselem} is again the key.  We argue by induction,
assuming that the Proposition holds for all $\Sigma_m$ sentences for all
$m\leq n$, uniformly in the sentence and the parameters.

Consider a positive sentence $(\exists\Fvec)~\psi(\Fvec,\avec)$,
where $\psi$ is $\Pi_n$ and positive.  The tree $T_{\exists\Fvec\psi,\avec}$ is an
$\avec$-decidable subtree of $\TQbar$, and we are given the parameters $\avec$
for our computation.  This tree
has as its paths exactly the realizations $\fvec$ of $\psi(\Fvec,\avec)$ in $\ab$,
so we simply use the Uniform Low Basis Theorem to compute an approximation
to the jump of a path through $T_{\exists\Fvec\psi,\avec}$.  This is a superapproximation
of a function $\G$ with $\ab\models\psi(\G(\avec),\avec)$.  Now the inductive
hypothesis, applied to the formula $\psi(\G(\avec),\avec)$ with the entire concatenated
tuple $\G(\avec)\hat{~}\avec$ viewed as the parameters, gives uniformly superapproximable generalized Skolem functions
for $\Sigma_{n-1}$ formulas, hence also gives a superapproximation of a Skolem function
for the $\Pi_n$ formula $\psi(\G(\avec),\avec)$, uniformly in $\avec$, and $\psi$.
Lemma \ref{lemma:super^2} completes the proof.
\comment{
For a positive $\Pi_{n+1}$ sentence $(\forall\Fvec)~\psi(\Fvec,\avec)$, the uniformity
in the inductive hypothesis is all that is needed.  Indeed, given $\psi$, $\avec$,
and any tuple $\fvec\in(\ab)^k$, we simply apply that hypothesis to the sentence
$\psi(\fvec,\avec)$, again regarding the entire string $\fvec\hat{~}\avec$ as the
tuple of parameters for this sentence.
}
\qed\end{pf}

\comment{
Our next theorem is stated for $\Pi_3$-separated formulas, but the argument
really applies only to $\Sigma_2$-separated formulas.  The parameters $\fvec$
in the initial string of $\forall$-quantifiers in the $\Pi_3$-separated formula $\theta(\avec)$
are simply folded in with the given parameters $\avec$.
}

Moving to non-positive formulas, we recall that, although Skolem functions
for $\Sigma_1$ sentences are merely constants, it must be proven that those constants
can be found uniformly in the sentence.

\begin{prop}
\label{prop:Sigma1}
$\Sigma_1$ formulas $\phi(\avec)$ have uniformly superapproximable generalized Skolem functions in $\ab$.
In other words, there is a procedure, uniform in $\phi$ and parameters $\avec$, that superapproximates a tuple
in $(\ab)^m$ of witnesses to the truth of each $\Sigma_1$ sentence $\phi(\avec)$ true in $\ab$.
\end{prop}
\begin{pf}
Say that $\phi(\avec)$ is of the form $(\exists\Gvec)[\vee_i\alpha_i(\Gvec,\avec)]$, where each
$\alpha_i$ is a conjunction of equations and inequations.  In Theorem \ref{thm:Sigma1},
we built an infinite subtree of $(\TQbar)^m$ and argued that it must contain a path of the
degree required there.  Here we have no specific degree requirement.  However,
the tree in Theorem \ref{thm:Sigma1} was constructed nonuniformly: it required
knowing some specific $i$ such that $\ab\models(\exists\Gvec)\alpha_i(\Gvec,\avec)$,
and also some specific level $n_0$ by which some witness in $\ab$ had satisfied
all of the inequations in that $\alpha_i$.  Each of these constituted finitely much information,
but here, those difficulties appear to contradict our uniformity claim.

Nevertheless, the theorem holds.  The reason is that a superapproximation
is able to change its mind finitely often and attempt to satisfy a different $\alpha_i$, 
or to use a different $n_0$, whereas a computation cannot do the same.
In the following uniform construction, which approximates a subtree
$T_{\phi,\avec}$ of $(\TQbar)^m$, there may exist stages $s$
(but only finitely many)
at which $T_{\phi,\avec,s}\neq T_{\phi,\avec,s+1}$.  The Uniform
Low Basis Theorem will give approximations $\Phi^{\avec}(e,s)$ at stage $s$
based on our stage-$s$ approximation of $T_{\phi,\avec}$.  This may produce
some strange values at early stages, but once we reach the stage $s_0$ after which
$T_{\phi,\avec,s}$ never changes again, it will behave as we desire from then on,
superapproximating the jump of a path through the limit tree $T_{\phi,\avec}$.
Thus the superapproximation will be uniform even though the construction of
the tree was not.

At each stage $s\geq 0$, we choose $\gamma_s$ to be the shortest node $\gamma\in(\TQbar)^m$
satisfying $\aut{F_l}\models\vee_i\alpha_i(\gamma,\avec\res F_l)$, where $l=\max(s,|\gamma|)$.
(If there are several such nodes at the same level, take the leftmost possibility as $\gamma_s$.)
Some such $\gamma_s$ must exist, since by hypothesis $\ab\models(\exists\Gvec)\vee_i\alpha_i(\Gvec,\avec)$
and for each witness $\gvec$, some finite initial segment of $\gvec$ will suffice.
Fix the least index $i_s$ such that $\aut{F_l}\models\alpha_{i_s}(\gamma,\avec\res F_l)$,
and for each $n$, define the $n$-th level of the tree $T_{\phi,\avec,s}$ to contain:
\begin{itemize}
\item
only $\gamma_s\res n$, if $n\leq |\gamma_s|$;
\item
those $\delta\in(\TQbar)^m$ of length $n$ such that
$$\gamma_s\sqsubseteq\delta~\&~
\an\models \alpha_{i_s}(\delta,\avec\res F_n),$$
if $n>|\gamma_s|$ and this set is nonempty; or
\item
all $\delta\in(\TQbar)^m$ of length $n$ such that
$\delta\res(n-1)\in T_{\phi,\avec,s}$, if $n>|\gamma_s|$ and the set
in the preceding item is empty.
\end{itemize}
Thus $T_{\phi,\avec,s}$ contains $\gamma_s$ but no other node at level $s$.
Above that level, we would like each level to contain all extensions of $\gamma_s$
that appear to witness satisfaction of the disjunct $\alpha_{i_s}$.  However,
in case we reach a level where no extension does so, we use the option
of extending the tree trivially instead.  This completes stage $s$.

Among all witnesses $\gvec\in(\ab)^m$ to the satisfaction of $\phi(\avec)$,
and all of the (finitely many) disjuncts $\alpha_i$ that each one satisfies,
there is some particular $\gvec$ and $i$ for which $\gvec$ satisfies $\alpha_i$
at the lowest level $l$.  (If several do this at the same level $l$,
fix some witness $\gvec$ with $\gvec\res F_l$ the leftmost possible.)
Then fix the least $i$ for which $\aut{F_l}\models \alpha_i(\gvec,\avec)$.
Eventually the construction will reach a stage $s_0$ with $\gamma_{s_0}=\gvec\res F_l$
and with this $i$ as $i_s$, and at all stages $s\geq s_0$ we will have
$\gamma_s=\gamma_{s_0}$ and $i_s=i_{s_0}$ and consequently
$T_{\phi,\avec,s}=T_{\phi,\avec,s_0}$.  Thus, as promised, our approximation
of this limit $T_{\phi,\avec}$ only changes at finitely many stages.

Now we apply the Uniform Low Basis Theorem to our approximation:
for each $e$ and $s$, let
$$\Phi^{\avec}(e,s) = \Upsilon^{T_{\phi,\avec,s}}(e,s).$$
By our construction, this superapproximates a path through $T_{\phi,\avec}$,
which is therefore a witness $\gvec$ with $\ab\models\vee_i\alpha_i(\gvec,\avec)$.
So we have the uniform generalized Skolem function required.
\qed\end{pf}

\begin{thm}
\label{thm:Skolem}
Every $\Sigma_2$-separated formula $\psi(\avec)$ has superapproximable
generalized Skolem functions, and the superapproximations
are uniform in the formula $\psi$.
\end{thm}
As described above, this proves that all $\Pi_3$-separated formulas
have superapproximable Skolem functions as well.
\begin{pf}
Here we refer to the construction of the tree $T_{\psi,\avec}$
in the proof of Theorem \ref{thm:compactlypos}.  Recall that this tree
is the limit of a sequence of subtrees $T_{\psi,\avec,s}$ of $(\TQbar)^m$.
At some stages these trees change,, but there was a finite
stage $s_0$ such that, for all $s\geq s_0$, $T_{\psi,\avec,s} = T_{\psi,\avec,s+1}$.
In general $s_0$ cannot be computed uniformly from $\psi$ and $\avec$.

In this theorem, $\psi(\avec)$ is of the form
$\exists\Gvec\lambda(\Gvec,\avec)$ with $\lambda$ $\Pi_1$-separated.
Assuming that $\ab\models\psi(\avec)$, we will produce
an $\avec$-computable superapproximation
$S_{\avec}(e,s)$ of a Skolem function $\G:(\ab)^k\to(\ab)^m$ for which 
$$ \ab\models \lambda(\G(\Fvec),\avec).$$
With $A=\oplus\avec$ as oracle,
the function $\G$ proceeds to compute the approximations
$T_{\psi,\avec,s}$ to the tree $T_{\psi,\avec}$,
just as in Theorem \ref{thm:compactlypos}.  It then applies the
Uniform Low Basis Theorem, so that for each $e$ at each stage $s$,
$S_{\avec}(e,s)$ superapproximates a path through the tree currently being built.
Specifically, at stage $s$, $S_{\avec}$ computes the greatest stage $s'\leq s$
at which the approximation to $T_{\psi,\avec,s'}$ reset itself
from the preceding stage instead of just extending the previous
$T_{\psi,\avec,s'-1}$.  $S_{\avec}$ treats stage $s'$ as the first stage
in the construction of $T_{\psi,\avec}$ and computes
its stage-$s$ approximation to a path through this tree using
the Uniform Low Basis Theorem.  Just as in Proposition \ref{prop:Sigma1},
the nonuniformity of the actual tree is simply folded into the process
of superapproximation.  Therefore,
$S_{\avec}$ is indeed a superapproximation of a path through $T_{\psi,\avec}$,
by which we mean that $\lim_s S_{\avec}(e,s)$ converges to the jump of such a path.
The path itself is the value $\G(\avec)$, which $S_{\avec}$ superapproximates.

All instructions here were entirely uniform in $\theta$ and in the parameters $\avec$,
so the uniformity claim in the theorem is clear.
\qed\end{pf}
\comment{
Notice that Theorem \ref{thm:Skolem} covers the special case of a $\Sigma_2$-separated
formula $\psi(\avec)$ as well, showing that one can compute witnesses $\gvec\in(\ab)^m$
to the truth of $\psi(\avec)$ uniformly in $\psi$ and $\avec$.  When one Skolemizes,
these witnesses are usually regarding as constants (or $0$-ary functions), not worthy
of computing, but initially it is not trivial to see that those constants can be found
uniformly in $\psi$ and $\avec$.  To see that they can be found uniformly, just run
the procedure in the theorem on the sentence $(\forall F)\psi(\avec)$, where the variable
$F$ does not appear in $\psi$. 
}

\comment{
Finally we claim that in $\ab$, due to its compactness (or more fully, to its presentation
as the paths through the finite-branching tree $\TQbar$), every sentence true in $\ab$
has Skolem functions that are Borel, indeed arithmetic.
}

\section{Questions}
\label{sec:questions}

The first question raised by our results, naturally, is how much further they can be extended.
Do the theorems of Sections \ref{sec:elementarity}, \ref{sec:definability}, and
\ref{sec:Skolem} extend to all $\Sigma_2$ formulas, or to $\Sigma_3$ formulas,
or beyond?  It does not seem possible to answer these questions with the techniques
employed here.  In current work, Jason Block is studying this situation for arbitrary
actions of compact groups on $\omega$ (as opposed to the specific case here,
where $\ab$ acts on the domain $\omega$ of $\Qbar$) and may be able to prove
broader results in that context.  Specific results for $\ab$ are likely to require
a stronger dose of field theory, unless the broader results turn out to resolve all questions.

On the other hand, one also naturally inquires into the subgroups $\abd$ of $\ab$,
defined by principal Turing ideals (containing all degrees $\leq\bfd$) rather than by
Scott ideals.  It was noted (as Theorem \ref{thm:KM} above) that Kundu and Miller have
shown that there is no computable generalized Skolem function even for $\Sigma_1$
formulas in $\ab$.  However, it remains possible (as of this writing) that every computable
$f\in\ab$ that happens to be of the form $g\circ g$ is in fact the square of a computable
$g\in \abz$:  all we know is that there is no uniform way of computing such a $g$
for each such $f$.  Thus, $\abz$ might yet prove to be an elementary subgroup
of $\ab$, at least for $\Sigma_1$ formulas and/or positive formulas and possibly more.
This will require real work, however.  At present we do not even know whether
$\ab$ and $\abz$ are elementarily equivalent!

Part of the interest in these issues revolves around the possibility of giving
a computability-theoretic description of $\ab$ as a direct limit.  Normally it is viewed
as a profinite group, of course, i.e., an inverse limit of finite groups, and this view
is totally natural and extemely productive.  Nevertheless, it could also be viewed
as the direct limit of the countable groups $\abd$, under inclusion, as $\bfd$
ranges over all Turing degrees.  It is not difficult to imagine that these two views
together might yield results that could not be derived from either one on its own.
This would seem especially likely if we can uncover further information about the
groups $\abd$.  How close are their theories to that of $\ab$?  How close
are they to being elementary?  To what extent does the isomorphism class
of $\abd$ depend on the choice of the degree $\bfd$?

This article has provided some evidence that the
Scott-ideal groups $\abI$ seem to be closer to elementary (within $\ab$)
than the principal-ideal groups $\abd$.  This being the case, it might
also be fruitful to consider $\ab$ as the direct limit of the groups $\abIA$,
as $A$ ranges over all subsets of $\omega$, or possibly over a class of representatives
of the Turing degrees, and $I_A$ is as in Theorem \ref{thm:exactpair}.
(That theorem also shows that $\abIA=\abb\cap\abc$, for the degrees
$\bfb$ and $\bfc$ of the two specific sets $B_A$ and $C_A$ defined there.
Thus each $\abIA$ has a simpler definition than the procedure in Subsection
\ref{subsec:Scottideals} would suggest.)
In order to make this direct limit close to effective,
various technical details are required, some of which appear in \cite{M23}.
Once again, the idea would be to take a direct limit of these subgroups
under simple inclusion $\subseteq$.  It is not clear, though, whether
one can build the subgroups $\abIA$ to respect Turing reducibility,
as it is unknown whether the $I_A$ built using our present version
of the Uniform Low Basis Theorem
satisfy $A\leq_T B \implies I_A\subseteq I_B$.  (This is expressed
in \cite{M23} as the question of whether the Uniform Low Basis Theorem
is or can be made \emph{monotone}.)

Of course, it has been known ever since the Lowenheim-Skolem
Theorems that there are countable elementary subgroups of $\ab$, and indeed that $\ab$
is the union of such subgroups.  The point of the statement here is the natural definition
of the subgroups:  by Theorem \ref{thm:exactpair}, each $\abIA$
can be defined by a simple statement of Turing reducibility
$f\leq_T B_A~\&~f\leq_T C_A$.  Additionally, one suspects that a group
such as $\abz$ has a presentation that is closer to computable than
any presentation of the natural L\"owenheim-Skolem subgroups would be.
Indeed, it seems likely that the usual L\"owenheim-Skolem procedure,
applied to produce a countable elementary subgroup of $\ab$ containing
a single given computable automorphism $f$ of $\Qbar$, may yield
the subgroup of all arithmetic automorphisms of $\Qbar$.
(Subgroups of $\ab$ defined by arithmetic reducibility, where
$$ A\leq_a B \iff A\text{~is arithmetic in~}B,$$
offer another potential presentation of $\ab$ as a direct limit,
this time most likely of elementary subgroups.  However, these
subgroups are much further from computable than the ones considered in this article.)

In the context of reverse mathematics, one naturally suspects that
all results here are at the level of $\textbf{WKL}_0$, due to the
frequent use of K\"onig's Lemma.  Dorais, Hirst, and Shafer
initiated exploration of such questions in \cite{DHS13}.
It would be surprising if any results here turned out to be stronger
than that, in terms of reverse mathematics, but we do not
address this question here.  For background in this area, see \cite{H15}.

\comment{

Could subgroups $\abd$ be elementary?

Describe the idea of $\ab$ as a direct limit.
Recalling that $\abd = \set{f\in\ab}{\deg f\leq \bfd}$ is the subgroup of all $\bfd$-computable
automorphisms of $\Qbar$, we notice that Theorem \ref{thm:exactpair} allows us to present
$\ab$ as an inverse limit of countable subgroups, each of which is elementary for compactly positive
sentences, as described in Section \ref{sec:elementarity}:
$$ \ab = \cup_{\bfa} (G_{\bfb_{\bfa}} \cap G_{\bfc_{\bfa}}),$$
where $\bfa=\deg{A}$ ranges over all Turing degrees (or as $A$ ranges over $2^\omega$)
and $\bfb_{\bfa}$ and $\bfc_{\bfa}$ are the degrees of the sets approximated uniformly in $A$
in Theorem \ref{thm:exactpair}.

CONJECTURE:  The underlying directed system is just the Turing degrees
under Turing reducibility.
%

}

\section{Appendix: computability and Scott ideals}
\label{sec:appendix}

\subsection{Scott ideals and the Low Basis Theorem}
\label{subsec:lowbasis}

A decidable infinite subtree of the complete binary tree
$2^{<\omega}$ need not have a computable path.
This follows from the Incompleteness Theorems, as the completions of a given
(consistent, decidable) axiom set form the set of paths through such a tree.
Thus the natural effective version of K\"onig's Lemma fails to hold.

More generally, a finite-branching infinite computable tree $T$ need not
have a computable path, and this remains true even if we assume that
the \emph{branching} in $T$ is computable, i.e., that we can compute
the function mapping each $\sigma\in T$ to the (finite) number of immediate
successors of $\sigma$ in $T$.  The above example with completions of axiom sets
already contradicts this stronger version.

With an oracle for the Halting Problem, one can readily compute paths
through infinite computable finite-branching trees.  However, in a surprise
for the logic community, this turned out not to be sharp.
Shoenfield showed in \cite{S60} that one can always do better,
and then Jockusch and Soare proved the following result in \cite{JS72}.
\begin{thm}[Low Basis Theorem, Jockusch and Soare]
\label{thm:lowbasis}
Every computable infinite tree with finite computable branching contains
a path $P$ of \emph{low} Turing degree, i.e., for which the jump
$P'$ is computable from the Halting Problem.
Indeed, there exists a single low degree $\bfd$ that can compute
a path through every such tree.  (Such a degree is now known
as a \emph{\PA-degree}).
\end{thm}
The term \emph{\PA-degree} arises from the fact that the
degrees with this property are precisely the Turing degrees of
complete extensions of the usual axiom set \PA\ for Peano Arithmetic.
However, the version above relativizes more readily, as follows.
\begin{defn}
\label{defn:PA}
For a given Turing degree $\bfc$ (or a set $C$ of this degree),
a Turing degree $\bfd$ is a \emph{\PA-degree relative to $\bfc$ (or to $C$)}
if $\bfd$ can compute a path through every $\bfc$-computable
infinite tree with $\bfc$-computable finite branching.
\end{defn}

Recall that a \emph{Turing ideal} is a proper nonempty subset of the Turing degrees
that is closed under the join operation and also closed downwards under
Turing reducibility.  This notion is often strengthened as follows.
\begin{defn}
\label{defn:Scottideal}
A \emph{Scott ideal} $I$ is a Turing ideal with the additional property
that, for every $\bfc\in I$, $I$ also contains a degree $\bfd$ that is
\PA\ relative to $\bfc$.
\end{defn}
Not all Turing ideals are Scott ideals.  Indeed, the \emph{principal}
Turing ideal $\set{\bfd}{\bfd\leq\bfc}$ generated by a single degree $\bfc$
can never be a Scott ideal, as it contains no \PA-degree relative to $\bfc$.
Scott ideals appear frequently in reverse mathematics:  every
$\omega$-model of the axiom system $\WKL$ defines a Scott ideal,
namely the set of degrees of functions $\omega\to\omega$
in the model.  Conversely every Scott ideal $I$ gives the
second-order part of an $\omega$-model of $\WKL$ in this way.
In practice, the virtue of a Scott ideal $I$ is that, whenever
$I$ contains both the degree of an infinite finite-branching tree $T$
and the degree of its branching, $I$ must also contain the degree
of some path through $T$.  In Section \ref{sec:elementarity}
we exploit this virtue.  For further discussion of the preceding concepts,
we recommend the book \cite{H15}.

\subsection{Producing Scott ideals}
\label{subsec:Scottideals}

We discuss how Scott ideals
may be found.  The natural intuition is to start with a degree $\bfa$ and build
an increasing sequence $\bfa=\bfa_0<\bfa_1<\cdots$ with each $\bfa_{j+1}$
being \PA\ relative to $\bfa_j$.  This succeeds, but of course it cannot be done effectively:
even $\bfa_1$ simply cannot be computed from $\bfa_0$.  Instead we appeal to the
\emph{Uniform Low Basis Theorem} of Brattka, de Brecht, and Pauly from \cite{BdBP12},
in the form given here.
\begin{thm}[Uniform Low Basis Theorem, in \cite{BdBP12}]
\label{thm:ULB}
There is a Turing functional $\Upsilon$ with the property that, for every
infinite subtree $T\subseteq 2^\omega$,
$\Upsilon^T:\omega^2\to\{ 0,1\}$ is total and there is a path $P$ through $T$ such that
$$(\forall n)~ \lim_s \Upsilon^T (n,s) = \chi_{P'}(n).$$
In the language of Definition \ref{defn:superapproximable},
$\Upsilon^T$ \emph{superapproximates} $P$.
\end{thm}
Here, since the jump $P'$ is the limit of a $T$-computable function, $P$ itself is certainly low
relative to $T$ (i.e., $P'\leq_T T'$).  Of course, $P$ itself is also the limit of a uniform
$T$-computable function (see \cite[Lemma 1]{M23} for a proof).  The true point, however,
is not that $P$ is uniformly $T$-limit-computable, but that $P'$ is:  this proves the lowness of $P$
relative to $T$, which was the point of the original (relativized but non-uniform) Low Basis Theorem
of Jockusch and Soare.

Theorem \ref{thm:ULB} allows us to produce the degree sequence $\{\bfa_j\}_{j\in\omega}$
as effectively as could be hoped.  (See \cite[Theorem 1]{M23} for details.)
Notice also that $\bfa'=\bfa_0' \geq \bfa_1'\geq\bfa_2'\geq\bfa_3'\geq\cdots$,
so every $\bfa_j$ (not just $\bfa_1$) is low relative to $\bfa$.  It now becomes possible to
apply Spector's notion of an \emph{exact pair}, and the construction there turns out to preserve
lowness relative to $\bfa$.  This result was proven much earlier using coded $\omega$-models
of \WKL, but the direct construction using the Uniform Low Basis Theorem,  which we quote here,
appears in \cite[Th.\ 4]{M23}.
\begin{thm}
\label{thm:exactpair}
For every subset $A\subseteq\omega$, there exist subsets $B$ and $C$ of $\omega$
whose join is low relative to $A$ (i.e., $(B\oplus C)'\leq_T A'$) and such that
$$\set{\deg{D}}{D\subseteq\omega\text{~with~}D\leq_TB~\&~D\leq_T C}$$
forms a Scott ideal $I_A$ containing $\deg{A}$.  (Thus every degree in this ideal
is low relative to $\deg A$.)  Furthermore, $A$-computable approximations
of such $B$ and $C$ may be constructed uniformly in the oracle $A$.
\end{thm}

The ideal $I_A$ is thus the intersection of the lower cones defined by the degrees of $B$ and $C$.
It is clear that $B$ and $C$ cannot have a meet, as that meet would have the greatest degree
in $I_A$, leaving nothing in the Scott ideal to be \PA\ relative to that degree.  However,
$\deg{A}$ comes reasonably close to being the greatest degree in $I_A$, as $B$ and $C$
(and even their join) are low relative to $A$.  For each specific set $A$,
we may sometimes denote these sets by $B_A$ and $C_A$ and their
degrees by $\bfb_A$ and $\bfc_A$.  We will refer to them as an
\emph{exact pair} for the ideal $I_A$, in keeping with the original terminology
of an exact pair for a countable increasing sequence of degrees.
(Indeed, the construction in \cite{M23} builds the sequence $\bfa_0<\bfa_1<\cdots$
described above and then takes $B$ and $C$ to be an exact pair for this sequence.)


\end{document}